\documentclass{amsart}

\usepackage{amsmath}
\usepackage{amsthm}
\usepackage{amssymb}

\usepackage{pict2e}

\newtheorem{thm}{Theorem}[section]
\newtheorem*{thm*}{Theorem}
\newtheorem{prop}[thm]{Proposition}

\newtheorem{defi}[thm]{Definition}
\newtheorem*{assum}{Assumption (A)}

\newcommand{\R}{\mathbb R}
\newcommand{\B}{\mathcal B}

\newcommand{\lap}{\mathcal L}
\newcommand{\schro}{\mathcal A}
\newcommand{\A}{\mathcal A}

\newcommand{\ie}{{\it i.e. }}

\begin{document}
\title[Interlacing property of Schr\"odinger operators]{Interlacing property of zeros of eigenvectors of Schr\"odinger operators on trees}
\author{Fran\c{c}ois Chapon}
\address{Universit\'e Paul Sabatier, Institut de Math\'ematiques de Toulouse, 118 route de Narbonne, F-31062 Toulouse Cedex 9}
\email{francois.chapon@math.univ-toulouse.fr}
\date{}
\begin{abstract}
We prove an analogue for trees of Courant's theorem on the interlacing property of zeros of eigenfunctions of a Schr\"{o}dinger operator. Let $\Gamma$ be a finite tree, and $\mathcal A$ a Schr\"{o}dinger operator on $\Gamma$. If the eigenvectors of $\mathcal A$ are  ordered according to increasing eigenvalues, and the vertices corresponding to zero coordinates are of degree at most two, then the zeros of the linear extensions of eigenvectors have the interlacing property.
\end{abstract}
\renewcommand{\subjclassname}{%
  \textup{2010} Mathematics Subject Classification}
\maketitle

\section{Introduction}

The famous Courant's theorem \cite{courant} about nodal domains of eigenfunctions of differential operators states that:
\begin{thm*}[Courant's theorem]\label{couthm}
Let $L$ be a self-adjoint second order elliptic operator on a domain $G$ with arbitrary boundary conditions. If its eigenfunctions are ordered according to increasing eigenvalues, then the nodes of the $n$-th eigenfunction $u_n$ divide the domain into no more than $n$ subdomains.
\end{thm*}
The subdomains defined in the theorem are called nodal domains and are the connected components  of the complement of the nodal set  $\{x\, | \, u_n(x)=0\}$, which are  separated by the nodes, or zeros, of eigenfunctions of $L$. 

If $L$ is a Schr\"odinger operator on $G=[a,b]\subset \R$, that is $L=\lap+V$ where $\lap$ is the Laplacian and $V$ some potential, Courant's theorem becomes: {\it The nodes of the $n$-th eigenfunction $u_n$ divide $G$ into exactly n nodal domains}. Besides, the zeros of eigenfunctions of $L$ have an interlacing property: {\it Between two zeros of $u_n$, there is exactly one zero of $u_{n+1}$.}

Analogues of Courant's result on graph have recently received  increasing attention, both in the mathematical and in the physical literature, and 
one may cite for instance the book \cite{biyibook} for a good introduction to this subject, and \cite{dav} for some historical comments.
In \cite{dav}, the authors Davies et al. prove a upper bound discrete nodal theorem  on graphs. 
Contrary to the situation when $G$ is a manifold, since a function on a graph is only defined on the set of vertices, an eigenfunction can change its sign without passing through zero, hence the nodes are not well defined.  So Davies et al. introduce the notion of strong sign graph, which is a connected set of vertices on which the eigenvector has the same sign, and prove the analogue of Courant's result on the maximal number of sign graphs of eigenvectors of a Schr\"odinger operator. Their proof is based on Courant's minimax theorem, see \cite{courant}, and some straightforward algebra.
In \cite{biyi}, B{\i}y{\i}ko\u{g}lu proves an analogue of Courant's theorem on trees under a certain genericity 
 condition:  
he proves that if $u$ is an eigenvector associated to the $n$-th ordered eigenvalue $\lambda_n$ and without a vanishing coordinate, then $u$ has exactly $n$ strong sign graphs 
and also gives a computational algorithm to find an eigenvector with minimum number of sign graphs when allowing multiplicity.
In \cite{berkolaiko}, Berkolaiko proves also a nodal theorem on graphs  under the same condition (both in the discrete and metric case). More precisely he gives a lower and upper bound of the nodal domains count of eigenfunctions of a generalized Laplacian on graphs and in particular recovers B{\i}y{\i}ko\u{g}lu's result on trees.
Recently, Xu and Yau \cite{xu-yau} give  a uniform proof of the lower and upper bounds of strong nodal domains on possibly disconnected graph, hence extending the results of \cite{berkolaiko}. 

On the physical level, the study of nodal domains appears in many different areas, such as quantum chaos, isospectral properties, or percolation theory. One may cite    for instance \cite{band-oren-smilansky} where several algorithmic and analytic  methods for counting nodal domains are presented,  \cite{gnutzmann-karageorge-smilansky} for the connection between  nodal domains count on manifolds  and geometrical content of the domain, 
\cite{oren-band-JPA}  for an example of isospectral graphs with the same nodal count sequence,
or \cite{smilansky-chaos} for nodal domains statistics as a criterion for quantum chaos.

Here, we are interested in the analogue of Courant's interlacing property  for Schr\"odinger operators on finite trees.  
Since on a tree there is a unique path connecting two vertices, we can extend a function by linearity on the edges, which allows us to define the zeros of an eigenfunction not only on vertices. This idea of geometric realization of a tree goes back to Friedman \cite{friedman}. To assure that the zeros of the eigenfunctions are well-defined, we make the assumption that the zeros coordinates of eigenvectors are of degree at most two, which implies that the eigenvalues are simple and that the corresponding eigenvectors are  without a zero graph.  The situation becomes then very similar to what happens in the real unidimensional case, 
and  we will prove in Theorem~\ref{thminterlacetree} the
interlacing property of the zeros 
 of the linear extension of eigenvectors of a Schr\"odinger operator. In particular we recover the exact nodal count already known for trees \cite{biyi,berkolaiko}.  The proof of these facts follows the lines of the classical proof of Courant \cite{courant} and are based on an analogue of 
  Green's formula
and on the discrete nodal upper bound proved in \cite{dav}.
 Note that, inspired by our method, Griffing, Lynch and Stone prove in \cite{griffing-lynch-stone-LAA} the interlacing property of the zeros   of harmonic functions of the Laplacian, defined as the extension of the eigenvectors of the Laplacian of a subgraph given by the Schur complement of the Laplacian with respect to the pendant vertices. 

The following is organized as follows. In section~\ref{Notations}, we introduce the notations, present the tree geometric realization allowing to consider linear extension of functions on edges, and discuss some examples. 
 Section~\ref{greensection} is devoted to prove the Green's formula on trees, and finally in section~\ref{interlacesection} we prove the Courant's nodal theorem on the interlacing property of zeros of eigenvectors of a Schr\"odinger operator on a tree.

\section{Definitions and notations}\label{Notations}

\subsection{Schr\"odinger operators}

Let $\Gamma=(V,E)$ be a finite tree, where $V$ is the set of vertices and $E$ the set of edges, with $|V|=N$, and $|E|=N-1$. Recall that a tree is a connected graph without cycles. This means that on a tree, any two points (possibly on the edges) are connected by a unique path. We note for $x,y\in V$, $x\sim y$ if there is an edge connecting $x$ and $y$, and $x$ and $y$ are said to be neighbors, or adjacent. The degree of a vertex $x$ is the number of edges connecting $x$. We choose some vertex 
 with only one neighbor and call it the root of the tree. This gives an orientation of the edges, the positive orientation being the direction connecting the root 
and the vertices. We note $(x,y)$ the edge starting from $x$ and ending at $y$. We consider weighted trees, that is there is a function $c\colon V\times V\rightarrow \R$ such that
\[
\begin{cases}
c(x,y)=c(y,x)>0,&\text{if $x\sim y$,}\\
c(x,y)=0,&\text{otherwise.}
\end{cases}
\]
For $x\sim y$, let $l(x,y)=\frac{1}{c(x,y)}$. As in \cite{friedman},  $l(x,y)$ is called the length of the edge connecting $x$ and $y$. 
Define $L^2(V)$ the space of functions on $V$ endowed with the scalar product
\[
\langle u,v \rangle_V=\sum_{x\in V}u(x)v(x),
\]
for $u,v$ functions on $V$. Define also $L^2(E)$ the space of functions $f$ on $E$ such that $f(x,y)=-f(y,x)$ for all edges $(x,y)\in E$, endowed with the scalar product
\[
\langle f,g \rangle_E=\sum_{e\in E} f(e)g(e)=\frac{1}{2}\sum_{x\in V}\sum_{y\in V,y\sim x}
f(x,y)g(x,y),
\]
for $f,g\in L^2(E)$, where we use the notation $f(x,y)=f\big((x,y)\big)$ for functions on $E$ for clarity.
The derivative operator $\partial\colon L^2(V)\rightarrow L^2(E)$ is defined by 
\[
\partial u(x,y)=c(x,y)^{1/2}(u(x)-u(y)),
\]
for $u\in L^2(V)$, and its adjoint 
$\partial^\ast\colon L^2(E)\rightarrow L^2(V)$ is then given by 
\[
\partial^\ast g(x)=\sum_{y\sim x}c(x,y)^{1/2} g(x,y), 
\]
for $g\in L^2(E)$.

Let $\lap=\partial^\ast\partial\colon L^2(V)\rightarrow L^2(V)$.
Then $\lap$ is called the {\it Laplacian} on $\Gamma$, and for $f\in L^2(V)$, we have for all $x\in V$,
\[
\lap f(x)=\sum_{y\in V\!\! ,\, y\sim x} c(x,y)\big(f(x)-f(y)\big). 
\]
Note that $\lap$ is a self-adjoint operator on $L^2(V)$.
If we see $f$ as a vector in $\R^N$, then $\lap$ can be seen as a $N\times N$ symmetric matrix whose nonzero entries are given by 
\[ 
\begin{cases} 
\lap_{xy}=-c(x,y),&\text{for $x\sim y$,}\\
\lap_{xx}=\sum_{y\sim x} c(x,y),&\text{on the diagonal.}
\end{cases}
\]
We  recall now the notion of Schr\"odinger operators. Let $r\colon V\rightarrow \R$ be some function on $V$, which plays the role of some  potential. Define $\schro=\lap +r$. Then, for all $f\in L^2(V)$, and all $x\in V$,
\[
\schro f(x)=\sum_{y\in V \!\! ,\, y\sim x}c(x,y)\big(f(x)-f(y)\big)+r(x)f(x).
\]
The operator $\schro$ is called a {\it Schr\"odinger operator}, or a generalized Laplacian, on $\Gamma$. As for the Laplacian, $\schro$ can be seen as a $N\times N$ symmetric matrix, with non-positive off-diagonal elements.

In the sequel, we will note $\lambda_i$, $i=1,\ldots,N$, the (real) eigenvalues of $\schro$ ordered such that 
\[
\lambda_1\leq\lambda_2\leq\cdots\leq\lambda_N,
\]
and such that the eigenspaces are orthogonal with respect to $\langle\cdot,\cdot \rangle_V$.
It is well known that, by the Perron-Frobenius theorem, $\lambda_1$ is simple and the first eigenvector can be chosen everywhere positive, see \cite{dav}. By orthogonality, any eigenvector associated with an eigenvalue different from $\lambda_1$ must then changes sign on $V$.

\subsection{Discrete nodal theorem}

We introduce  now the notion of sign graphs of a function on $V$ as defined in \cite{dav}. 
\begin{defi}
Let $u$ be a function on $V$.
A strong positive (resp. negative) sign graph of $u$ is a maximal subtree $S$ of $\Gamma$ with $u(x)>0$ (resp. $u(x)<0$), for all vertices $x$ of $S$.
\end{defi}

Then the theorem of Davies et al. \cite{dav}, which is an analogue on graph of the Courant's nodal theorem, states that:
\begin{thm}[\cite{dav}]
Let $\lambda_1,\ldots,\lambda_N$ be the ordered eigenvalues of a Schr\"odinger operator on $\Gamma$. Suppose $\lambda_n$ is of multiplicity $r$. Then any eigenvector corresponding to $\lambda_n$ has at most $n+r-1$ strong sign graphs.
\end{thm}

Let $u$ be an eigenvector of $\schro$ with eigenvalue  $\lambda$. Let $x$ be a {\it zero vertex} of $u$, \ie $u(x)=0$. Then, $\schro u(x)=0$, and we have
\[
\sum_{y\sim x}c(x,y)u(y)=0.
\]
Since $c(x,y)>0$ for all $x\sim y$, two cases are possible: we have either $u(y)=0$ for all $y\sim x$, and $x$ is said to belong to a {\it zero graph}, or $x$ is adjacent to both strict signs.
For our context of Schr\"odinger operators, we make the following assumption.
\begin{assum} \label{assumpnodal}
For all $i=1,\ldots,N$, denote by $u_i$ an eigenvector  associated with $\lambda_i$.   Let $s_k^{(i)}$ be the number of vertices corresponding to  the zeros of $u_i$ and of degree $k$. Then $s_k^{(i)}=0$ for all $k\geq3$ and all $i=1,\ldots,N$, that is any zero vertex is of degree at most two.
\end{assum}
This assumption implies that if there is a zero vertex, it is of degree exactly two. Indeed, if there is a zero graph, then since a tree is connected it has to end with a zero vertex adjacent to both strict signs, thus of degree at least  three, hence contradicting the assumption. A zero vertex has then exactly two adjacent vertices of both strict signs and no end vertex can be a zero vertex.
By a result of Fiedler, see Theorem~2.4 and Corollary~2.7 in \cite{fiedler}, this is equivalent of saying that all eigenvalues are simple and the corresponding eigenvectors are without a zero graph.  This will insure that the zeros of the linear extension of the eigenvectors are well defined, see next subsection. 

Under assumption~(A), the theorem of Davies et al. becomes then: any eigenvector associated  with $\lambda_n$ has at most $n$ strong sign graphs.  Note that B{\i}y{\i}ko\u{g}lu \cite{biyi} (see also \cite{berkolaiko}) proves under the  assumption that  all eigenvectors are without  vanishing coordinates, that any eigenvector associated with $\lambda_n$ has exactly $n$ strong sign graphs.  %



\subsection{Geometric realization}

As introduced by Friedman \cite{friedman} (see also \cite{biyibook}), we consider a geometric realization of $\Gamma$. An edge $(x,y)$ can be identified with the real interval $[0,l(x,y)]$.
A function $u$ defined on $V$ can then be extended on the edges by linearity. Let us call $\tilde u_{xy}$ this extension on the edge $(x,y)$. We have then
\[
\tilde u_{xy}(t)=
\frac{u(y)-u(x)}{l(x,y)}t+u(x),\qquad\text{for $t\in [0,l(x,y)]$.}
\]
We then denote by $\tilde u$ the linear extension of $u$ on each edge.


We introduce now the correspondence between strong sign graphs of eigenvectors of $\schro$ and the nodal domains of their respective linear extensions.

\begin{defi}
Let $u$ be a function on $V$, and $\tilde u$ its linear extension on edges.
A nodal domain of $\tilde u$ is a maximal connected path on which $\tilde u$ does not vanish. 
\end{defi}
Since  a tree is connected,  this notion is well defined. For $u$ an eigenvector of $\schro$, the different nodal domains of $\tilde u$ are then separated by the zeros of $\tilde u$ which are well defined by assumption~(A) and possibly on the edges. Note that
a nodal domain of $\tilde u$ can also be seen as a subtree of $\Gamma$ with some incomplete edges.

Let $G$ be a strong sign graph of $u$. Let $(x,y)$ be an edge with $x\in G$ and $y\in V\setminus G$.
Then $u(x)u(y)\leq0$ on $(x,y)$, otherwise $y$ would be in $G$. This implies that $\tilde u$ must  vanish on the interval $]0,l(x,y)]$.  Hence, $G$ defines uniquely a nodal domain $\tilde G$ of $\tilde u$, because, since $\tilde u$ is linear on any edge, we cannot have two zeros on the same edge, and a nodal domain contains at least one vertex. So, there is a one-to-one correspondence between strong sign graphs of $u$ and nodal domains of its linear extension $\tilde u$.  

\subsection{Some examples}

\begin{figure}[h]
\begin{center}
\setlength{\unitlength}{0.8cm}
\begin{picture}(9,1)(0,-0.5)
\put(0,0){\line(1,0){2}}
\put(0,0){\circle*{0.2}}
\put(1,0){\circle*{0.2}}
\put(2,0){\circle*{0.2}}
\put(-0.14,0.25){\footnotesize $+$}
\put(0.86,0.25){\footnotesize $+$}
\put(1.86,0.25){\footnotesize $+$}
\put(0.86,-0.7){$u_1$}
\put(4,0){\line(1,0){2}}
\put(4,0){\circle*{0.2}}
\put(5,0){\circle*{0.2}}
\put(6,0){\circle*{0.2}}
\put(3.86,0.25){\footnotesize $-$}
\put(4.9,0.25){\footnotesize $0$}
\put(5.86,0.25){\footnotesize $+$}
\put(4.86,-0.7){$u_2$}
\put(8,0){\line(1,0){2}}
\put(8,0){\circle*{0.2}}
\put(9,0){\circle*{0.2}}
\put(10,0){\circle*{0.2}}
\put(7.86,0.25){\footnotesize $+$}
\put(8.85,0.25){\footnotesize $-$}
\put(9.86,0.25){\footnotesize $+$}
\put(8.86,-0.7){$u_3$}
\end{picture}
\caption{Eigenvectors of a linear graph with weights $1,1$.}\label{eigenlineartree}
\end{center}
\end{figure}
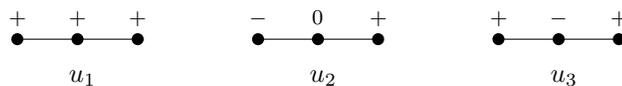

First let us consider the simple example of a linear graph with weights $1,1$ on both edges. 
The Laplacian is given by
\[
\mathcal L =\begin{pmatrix} 1 & -1 & 0 \\ -1 & 2 & -1 \\ 0 & -1 & 1  \end{pmatrix}.
\]
The eigenvalues are $0$, $1$ and $3$ and are simple (recall that it is a well known fact that the eigenvalues of tridiagonal symmetric matrices with non-zero subdiagonal elements are simple), with eigenvectors given by 
\[
u_1=\begin{pmatrix} 1 \\ 1 \\ 1 \end{pmatrix}, \ \  u_2=\begin{pmatrix} -1 \\ 0 \\1 \end{pmatrix}, \ \ u_3=\begin{pmatrix}   1 \\ -2 \\ 1 \end{pmatrix},
\]
represented in Figure~\ref{eigenlineartree}.
 Thus we see that in the only nodal domain of $\tilde u_1$ there is exactly one zero of $\tilde u_2$, and in each nodal domain of $\tilde u_2$ there is exactly one zero of $\tilde u_3$ (both of them being on  an edge).

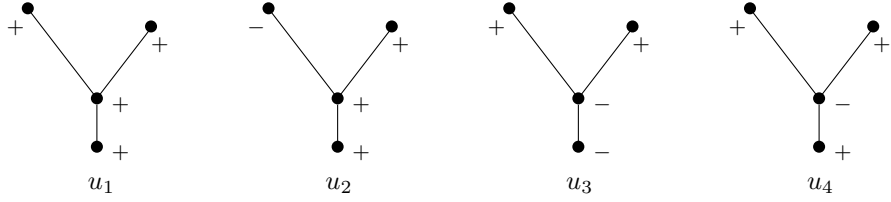
\begin{figure}[t]
\begin{center}
\setlength{\unitlength}{0.8cm}
\begin{picture}(12,3.5)(0,0)
\multiput(0,0.7)(4,0){4}{\line(0,1){0.7}}
\multiput(0,1.5)(4,0){4}{\line(1,1.3){0.9}}
\multiput(0,1.5)(4,0){4}{\line(-1,1.3){1.2}}
\multiput(0,0.7)(4,0){4}{\circle*{0.2}}
\multiput(0,1.5)(4,0){4}{\circle*{0.2}}
\multiput(0.9,2.7)(4,0){4}{\circle*{0.2}}
\multiput(-1.15,3)(4,0){4}{\circle*{0.2}}
\put(0.25,0.5){\footnotesize $+$}
\put(0.25,1.3){\footnotesize $+$}
\put(0.9,2.3){\footnotesize $+$}
\put(-1.5,2.6){\footnotesize $+$}
\put(4.25,0.5){\footnotesize $+$}
\put(4.25,1.3){\footnotesize $+$}
\put(4.9,2.3){\footnotesize $+$}
\put(2.5,2.6){\footnotesize $-$}
\put(8.25,0.5){\footnotesize $-$}
\put(8.25,1.3){\footnotesize $-$}
\put(8.9,2.3){\footnotesize $+$}
\put(6.5,2.6){\footnotesize $+$}
\put(12.25,0.5){\footnotesize $+$}
\put(12.25,1.3){\footnotesize $-$}
\put(12.9,2.3){\footnotesize $+$}
\put(10.5,2.6){\footnotesize $+$}
\put(-0.15,0){$u_1$}
\put(3.8,0){$u_2$}
\put(7.8,0){$u_3$}
\put(11.8,0){$u_4$}
\end{picture}
\caption{Eigenvectors of a star graph with weights $1,\frac14,\frac13$.}\label{starOK}
\end{center}
\end{figure}

Consider now the tree which is a star graph with four vertices and weights $1,\frac14,\frac13$. The associated Laplacian is given by
\[
\mathcal L = \begin{pmatrix}
1 & -1 & 0 & 0 \\ 
-1& \frac{19}{12} & -\frac14 & -\frac13 \\
0&  -\frac14& \frac14& 0 \\
0 &-\frac13 &0& \frac13
\end{pmatrix},
\]
and the eigenvalues are $0$, $\frac43 -\frac{\sqrt{10}}{3}$, $\frac12$, $\frac43 +\frac{\sqrt{10}}{3}$ and are simple. The eigenvectors are represented in Figure~\ref{starOK} and can be seen to have no zero vertices. Again one can see the interlacing property, that is in each nodal domain of $\tilde u_i$ there is exactly one zero of $\tilde u_{i+1}$.

\begin{figure}[h]
\begin{center}
\setlength{\unitlength}{0.8cm}
\begin{picture}(12,3.2)(0,-0.8)
\multiput(0,0)(4,0){4}{\line(0,1){1.4}}
\multiput(0,1.5)(4,0){4}{\line(1,1){0.6}}
\multiput(0,1.5)(4,0){4}{\line(-1,1){0.6}}
\multiput(0,0)(4,0){4}{\circle*{0.2}}
\multiput(0,1.5)(4,0){4}{\circle*{0.2}}
\multiput(0.6,2.1)(4,0){4}{\circle*{0.2}}
\multiput(-0.6,2.1)(4,0){4}{\circle*{0.2}}
\put(0.25,-0.2){\footnotesize $+$}
\put(0.25,1.25){\footnotesize $+$}
\put(0.8,2){\footnotesize $+$}
\put(-1.1,2){\footnotesize $+$}
\put(4.25,-0.2){\footnotesize $-$}
\put(4.25,1.25){\footnotesize $+$}
\put(4.8,2){\footnotesize $+$}
\put(2.9,2){\footnotesize $+$}
\put(8.25,-0.2){\footnotesize $0$}
\put(8.25,1.25){\footnotesize $0$}
\put(8.8,2){\footnotesize $+$}
\put(6.9,2){\footnotesize $-$}
\put(12.25,-0.2){\footnotesize $+$}
\put(12.25,1.25){\footnotesize $-$}
\put(12.8,2){\footnotesize $+$}
\put(10.9,2){\footnotesize $+$}
\put(-0.15,-0.7){$u_1$}
\put(3.8,-0.7){$u_2$}
\put(7.8,-0.7){$u_3$}
\put(11.8,-0.7){$u_4$}
\end{picture}
\caption{The star graph with weights $\frac12,1,1$ counterexample: $u_3$ has a zero graph and only two strong sign graphs.}\label{star3counterex}
\end{center}
\end{figure}
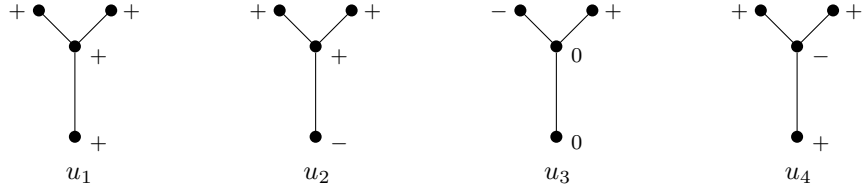

To show that assumption~(A) is essential, we provide the following counterexample, represented in Figure \ref{star3counterex}. 
Consider the star graph with four vertices, and weights $\frac12,1,1$. The ordered eigenvalues of the associated Laplacian are $0,2-\sqrt2,1,2+\sqrt2$, and in particular they are simple. The corresponding eigenvector for the third eigenvalue 
$1$  is $u_3=\vphantom{A}^t(0,0,-1,1)$ and has a zero graph, hence its linear extension is identically zero on one of the nodal domains of $\tilde u_2$, which fails Theorem~\ref{thminterlacetree}. 
Note also that $u_3$  has only 2 ($<3$) strong sign graphs. 

Counterexamples of eigenvectors associated with non-simple eigenvalue $\lambda_n$ having more than $n$ strong sign graphs can be found in \cite{dav}. 
For instance, one can show that the second eigenvalue of the Laplacian associated with a star graph with $N$ vertices with weight $1$ on each edge has multiplicity $N-2$, with an eigenvector having  only a zero vertex at the unique internal vertex (of degree $N-1$), and hence having  $N-1$ strong sign graphs, see \cite{dav} for more details.

\section{Green's formula}\label{greensection}

Analogues of Green's formula on graphs have been already considered, see for instance  \cite{biyibook}. To be self-contained we give a straightforward proof in our  context.  

Let us first introduce some boundary sets.
Let $u$ be some function on $V$. 
Let $G$ be a strong positive sign graph of $u$, and $\tilde G$ the corresponding  nodal domain  of $\tilde u$. Define the following boundary sets:
\begin{align*}
\delta(G)&=\left\{x\in V\setminus G \, | \, x\sim y, \text{ for some $y\in G$}\right\},\\
\partial(G)&=\left\{(x,y)\in E \, |\,  x\in G, y\in\delta(G)\right\}.
\end{align*}
Since $u(x)>0$ for all $x\in G$, we have $u(y)\leq 0$ for all $y\in \delta(G)$, otherwise $y$ would be in $G$.
This implies that $\tilde u_{xy}$ vanishes on $]0,l(x,y)]$, where $(x,y)$ is an  edge belonging to $\partial(G)$.
Let us call $t(x,y)$ the point of $]0,l(x,y)]$ where $\tilde u_{xy}$ vanishes. We have 
\[
t(x,y)=\frac{l(x,y)u(x)}{u(x)-u(y)}.
\]
The boundary of $\tilde G$  (depending on $u$) is  then defined by
\[
\B(\tilde G)=\left\{t(x,y)\, |\, (x,y) \in \partial(G)\right\}.
\]
This boundary set is defined in the same way for $G$ a strong negative sign graph.
Let $\nabla$ be the usual gradient on $\R$. We have the analogue of the Green's formula:

\begin{prop}[Green's formula]
Let $\schro=\lap +r$ be a Schr\"odinger operator.  Let $u,v$ be functions on $V$, and $\tilde u,\tilde v$ their linear extensions on edges. Let $G$ be a strong sign graph of $u$, and $\B(\tilde G)$ the boundary of the corresponding nodal domain of $\tilde u$. Then, we have,
\[
\sum_{x\in G}\schro u(x)v(x)-\sum_{x\in G}u(x)\schro v(x)=-\sum_{t\in \B(\tilde G)}\nabla\tilde u(t)\tilde v(t).
\]
\end{prop}
\begin{proof}
Since,
\[
\sum_{x\in G}\schro u (x)v(x)-\sum_{x\in G}u(x)\schro v(x)=\sum_{x\in G}\lap u(x)v(x)-\sum_{x\in G}u(x)\lap v(x),
\]
it suffices to prove the proposition for the Laplacian $\lap$.
 
We have,
\begin{align*}
\sum_{x\in G}\lap u(x)v(x)
&=
\sum_{x\in G}\sum_{\substack{y\in G\\ y\sim x}}c(x,y)\big(u(x)-u(y)\big)v(x)\\
&\qquad\qquad
+\sum_{x\in G}\sum_{\substack{y\in \delta(G)\\ y\sim x}}c(x,y)\big(u(x)-u(y)\big)v(x).
\end{align*}
But, 
\begin{align*}
\sum_{x\in G}\sum_{\substack{y\in G\\ y\sim x}}&c(x,y)\big(u(x)-u(y)\big)v(x)\\
&=\frac12\sum_{x\in G}\sum_{\substack{y\in G \\ y\sim x}}c(x,y)\big(u(x)-u(y)\big)\big(v(x)-v(y)\big),
\end{align*}
and the same holds by exchanging the roles of $u$ and $v$. Thus, we obtain
\begin{align*}
\sum_{x\in G}\lap u(x)v(x)-\sum_{x\in G}\lap v(x)u(x)
&= \sum_{x\in G}\sum_{\substack{y\in \delta(G) \\ y\sim x}}c(x,y)\big(u(x)-u(y)\big)v(x) \\
&\qquad\qquad-\sum_{x\in G}\sum_{\substack{y\in \delta(G) \\ y\sim x}}c(x,y)\big(v(x)-v(y)\big)u(x).
\end{align*}
On the interval $[0,l(x,y)]$, we have $\tilde v(t)=\frac{v(y)-v(x)}{l(x,y)}t+v(x)$, hence 
\[
\tilde v(t(x,y))=-\frac{v(x)-v(y)}{u(x)-u(y)}u(x)+v(x),
\]
for $t(x,y)=\frac{l(x,y)u(x)}{u(x)-u(y)}\in \mathcal B(\tilde G)$. Furthermore, since $\nabla$ is the usual derivate, we have 
\[
\nabla\tilde u(t)=\frac{u(y)-u(x)}{l(x,y)}=-c(x,y)\big(u(x)-u(y)\big), 
\]
for all $t\in [0,l(x,y)]$.
Hence, we have
\begin{align*}
\sum_{t\in\B(\tilde G)}\nabla\tilde u(t)\tilde v(t)
&=\sum_{x\in G}\sum_{\substack{y\in\delta(G)\\ y\sim x}}\nabla\tilde u(t(x,y))\tilde v(t(x,y))\\
&=\sum_{x\in G}\sum_{\substack{y\in\delta(G)\\ y\sim x}}  c(x,y)\big(v(x)-v(y)\big)u(x)\\
&\qquad\qquad-\sum_{x\in G}\sum_{\substack{y\in\delta(G)\\ y\sim x}} c(x,y)\big(u(x)-u(y)\big)v(x),
\end{align*}
so the proposition is proved.
\end{proof}

\section{Interlacing property of zeros of eigenvectors}\label{interlacesection}

As in the classical result  of Courant, we can now prove  the interlacing property of the zeros of eigenvectors of $\schro$.

\begin{thm} \label{thminterlacetree}
Let $\Gamma$ be a finite tree with $N$ vertices, and
 $\A$  a Schr\"odinger operator on $\Gamma$.
Let $\lambda_n,n=1,\ldots,N$, be the ordered eigenvalues of $\schro$, and let $u_n$ be an eigenvector associated with $\lambda_n$, and denote by $\tilde u_n$ its linear extension. If the zero vertices of the eigenvectors $u_n$'s  are of degree at most two, then the zeros of the  $\tilde u_n$'s interlace, in the sense that in any nodal domain of $\tilde u_{n-1}$ there is exactly one zero of $\tilde u_{n}$.
\end{thm}
\begin{proof}
Recall that we have seen that by Fiedler \cite{fiedler} if the zero vertices of the eigenvectors of $\A$ are  of degree at most two, then the eigenvalues of $\A$ are simple and the corresponding eigenvectors are without a zero graph.

Let $\lambda$ and $\mu$ be two eigenvalues of $\schro$ with $\lambda<\mu$, and let $u$ and $v$ be the associated eigenvectors. 
Let $G$ be a strong positive sign graph of $u$, and $\tilde G$ be the corresponding nodal domain of $\tilde u$, with boundary $\B(\tilde G)$.
Suppose that $\tilde v$ does not change sign on $\tilde G$, and say $\tilde v>0$ on $\tilde G$. This implies that $v(x)>0$, for all $x\in G$.  By Green's formula, we have
\[
(\lambda-\mu)\sum_{x\in G}u(x)v(x)=-\sum_{t\in\B(\tilde G)}\nabla\tilde u(t)\tilde v(t).
\]
The left hand side of the above expression is then negative. Since $\nabla\tilde u=c(x,y)(u(y)-u(x))< 0$ on the edges $(x,y)\in\partial G$, and $\tilde v>0$ on $\tilde G$, then the right hand side is non-negative. So there is a contradiction, and $\tilde v$ must vanish  and thus change sign on $\tilde G$ since a zero of $\tilde v$ is between two opposite sign vertices, even if the zero is a vertex since we cannot have zero graphs by assumption. Using the same argument on all of the sign graphs of $u$, we conclude that $\tilde v$ has at least one more nodal domain than $\tilde u$. By the discrete nodal theorem of \cite{dav}, since all eigenvalues $\lambda_i$ of $\mathcal A$ are simple, then $u_n$ has at most $n$ strong sign graphs, and hence $\tilde u_n$ has at most $n$ nodal domains. 
So, we deduce from this by iteration, that $\tilde u_n$ has exactly $n$ nodal domains, and hence $u_n$ has exactly $n$ strong sign graphs.

Since on a tree, there is a unique path connecting any two points of the tree, the different nodal domains of $\tilde u_n$ are separated by a unique zero of $\tilde u_n$. 
Since $\tilde u_n$ has exactly $n$ nodal domains separated by its zeros, $\tilde u_n$ has exactly $n-1$ zeros in $\Gamma$. 
Hence, since $\tilde u_n$ must vanish in the interior of every nodal domain of $\tilde u_{n-1}$,
the zeros of the  $(\tilde u_n)_{1\leq n\leq N}$ must interlace in the sense that  in any nodal domain of $\tilde u_{n-1}$, there is exactly one zero of $\tilde u_{n}$.
\end{proof}

\subsection*{Acknowledgments}
 
I would like to thank Gregory Berkolaiko for some useful comments and suggestions.


\providecommand{\bysame}{\leavevmode\hbox to3em{\hrulefill}\thinspace}
\providecommand{\MR}{\relax\ifhmode\unskip\space\fi MR }
\providecommand{\MRhref}[2]{%
  \href{http://www.ams.org/mathscinet-getitem?mr=#1}{#2}
}
\providecommand{\href}[2]{#2}

\end{document}